\documentclass[a4paper,12pt]{article}

\usepackage{graphicx}%
\usepackage{multirow}%
\usepackage{amsmath,amssymb,amsfonts}%
\usepackage{amsthm}%
\usepackage{mathrsfs}%
\usepackage[title]{appendix}%
\usepackage{xcolor}%
\usepackage{textcomp}%
\usepackage{manyfoot}%
\usepackage{booktabs}%
\usepackage{algorithm}%
\usepackage{algorithmicx}%
\usepackage{algpseudocode}%
\usepackage{listings}%

\usepackage{tikz}
\usetikzlibrary{arrows,decorations.markings,patterns}
\usetikzlibrary{arrows.meta,patterns.meta}

\newtheorem{theorem}{Theorem}[section]
\newtheorem{lemma}[theorem]{Lemma}
\newtheorem{cor}{Corollary}[theorem]

\begin{document}
\title{On the minimum number of entries in a set of maximal mutually orthogonal partial Latin squares}

\author{Diane M. Donovan\footnote{Donovan acknowledges the support of the Australian Government through funding of the Australian Research Council  Centre of Excellence for Plant Success in Nature \& Agriculture (project number CE200100015).},\\
\normalsize{ARC Centre of Excellence for Plant Success in Nature and Agriculture,}\\
\normalsize{School of Mathematics and Physics, University of Queensland,}\\
\normalsize{Brisbane 4072, Australia.} \\
\texttt{(dmd@maths.uq.edu.au)}\\
\\
Mike J. Grannell,\\
\normalsize{School of Mathematics and Statistics, The Open University,} \\
\normalsize{Walton Hall, Milton Keynes MK7 6AA, United Kingdom.}\\
\texttt{(m.j.grannell@open.ac.uk)}\\
\\
Emine \c{S}ule Yaz{\i}c{\i}\footnote{Yaz{\i}c{\i} acknowledges the support of the Turkish Government through funding by The Scientific and Technological Research Council of T\"{u}rkiye (T\"{U}B\.{I}TAK grant number: 121F111)},\\
\normalsize{Department of Mathematics, Ko\c{c} University,}\\
\normalsize{Sar{\i}yer, 34450, \.{I}stanbul, T\"{u}rkiye.}\\
\normalsize{Also}\\
\normalsize{School of Mathematics and Physics, University of Queensland,}\\
\normalsize{Brisbane 4072, Australia.} \\
\texttt{(eyazici@ku.edu.tr)}\\
}
\maketitle

\newpage



\begin{abstract}
In superimposed form, a set of $k\ge2$ mutually orthogonal partial Latin squares of order $n$ is said to be maximal if no additional $k$-tuple may be inserted into any empty cell without violating orthogonality or the Latin property. For $k=2$ and $n\ge 21$ it is shown that the minimum number of filled cells in such a square is $\lceil n^2/3 \rceil$. Likewise for $k=3$ and $n\ge 44$ the minimum number of filled cells is $\lceil n^2/4 \rceil$ or $\lceil n^2/4 \rceil+1$, depending on the residue of $n$ modulo 4. In each case it is shown that the extremal systems are essentially those produced by a construction presented in the paper. A similar result is obtained for all $k\ge4$ for minimal systems that have a certain balance property, and it is conjectured that for sufficiently large $n$, even without this balance property, the extremal systems always correspond to those produced by this construction.
\end{abstract}

\section{Introduction} \label{sec:intro}
The problem addressed in this paper concerns sets of $k$ mutually orthogonal partial Latin squares of order $n$. There is a connection to a Tur\'{a}n-type problem in graph theory, which is explained in more detail later in this Introduction. The general conjecture is that, for sufficiently large $n$, the minimum number of filled cells in a maximal set of $k$ mutually orthogonal partial Latin squares of order $n$ is approximately $n^2/(k+1)$, with the precise number depending on the residue of $n$ modulo $k+1$. The result for $k=1$ (a single square) was resolved by Horak and Rosa in 1993 \cite{HorakRosa}. Our 2025 paper \cite{MOPLS1} dealt mainly, but incompletely, with the case $k=2$. In the current paper, for $n\ge 21$, we completely resolve the case $k=2$ and identify the extremal systems. We then provide a similarly complete solution for $k=3$ when $n\ge 44$. Further progress is also made when $k\ge 4$.

We start by recalling some basic definitions. A \emph{Latin square} of order $n$, denoted by LS($n$), is an $n\times n$ array $L=[L(i,j)]$, with entries from an $n$-element set, such that each element of this set occurs once in every row and once in every column of $L$. If convenient, rows and columns can be indexed by $[n]=\{1,2,\ldots,n\}$, and the set of entries can also be taken as $[n]$.

A \emph{partial} Latin square of order $n$, denoted by PLS($n$), is defined in the same manner as a Latin square of order $n$ except that some of the cells may be empty, in other words, each element of the entry set occurs \emph{at most} once in every row and \emph{at most} once in every column.

A \emph{maximal} PLS($n$) is a partial Latin square of order $n$ that cannot be extended to another PLS($n$) by inserting any entry into any empty cell. A maximal PLS($n$) is denoted by MPLS($n$). Horak and Rosa \cite{HorakRosa} proved that a PLS($n$) with less than $n^2/2$ filled cells cannot be maximal.

For $k\ge 2$, $k$ Latin squares of the same order $n$ are said to be \emph{mutually orthogonal} if, when superimposed, each of the $\binom{k}{2} n^2$ possible ordered pairs of entries appears once and only once in the cells of the superimposed square. We denote such an array as a $k$-OLS($n$) ($k$-Orthogonal Latin Square of order $n$). Thus a $k$-OLS($n$) is a single array with each cell having an ordered $k$-tuple of entries. When a 2-OLS($n$) is presented as two separate squares, the squares are usually referred to as a pair of MOLS($n$) (Mutually Orthogonal Latin Squares of order $n$).

For $k\ge 2$, $k$ partial Latin squares of the same order $n$ and having the same filled cells as one another, are said to be \emph{mutually orthogonal} if, when superimposed, each of the $\binom{k}{2} n^2$ possible ordered pair of entries appears \emph{at most} once in the cells of the superimposed square. We denote such an array as a $k$-OPLS($n$) ($k$-Orthogonal Partial Latin Square of order $n$).

For $k\ge2$ a \emph{maximal} $k$-OPLS($n$) is a $k$-OPLS($n$) that cannot be extended to another $k$-OPLS($n$) by inserting any possible entry $k$-tuple into any empty cell. We denote a maximal $k$-OPLS($n$) as a $k$-MOPLS($n$). In our previous paper \cite{MOPLS1} when $k=2$, the ``2'' was dropped and the array was simply denoted as an MOPLS($n$). In that paper the focus was almost entirely on the case $k=2$. In the current paper we will retain the ``2'' to avoid confusion because we will obtain results also for $k=3$ and for general $k$.

With the exceptions of $n=2$ and $n=6$, the maximum size (number of filled cells) of a 2-MOPLS($n$) is $n^2$ since there is a pair of orthogonal Latin squares of order $n$ for every positive integer $n\ne2,6$. More generally, for fixed $k>2$ and all sufficiently large $n$, the maximum size of a $k$-MOPLS($n$) is also $n^2$ since there are $k$ mutually orthogonal Latin squares of order $n$ if $n$ is sufficiently large (see \cite[p.163]{HBK}). The focus in this paper is on the opposite extreme: seeking the minimum number of filled cells in a $k$-MOPLS($n$). Results concerning the structure of $k$-MOPLS($n$) may prove applicable to the problem of embedding a $k$-OPLS($n$) into a $k$-OLS($n'$) for some $n'\ge n$ (see \cite{DY}).

Provided that $n$ is sufficiently large, a $k$-MOPLS($n$) may be constructed as follows. The construction also works for $k=1$ and then gives an MPLS($n$).

\medskip

\noindent\textbf{Construction 1.}~
Suppose that $n=(k+1)s+r$, where $0\le r \le k$. Starting from an empty $n\times n$ array, place $(k+1-r)$ ~$k$-OLS($s$) arrays and $r$ \linebreak $k$-OLS($s+1$) arrays (all on different entry sets) down the leading diagonal. The ingredients exist provided that $n$ (and hence $s$) is sufficiently large. The number of filled cells $F$ is given by
\[F=(k+1-r)s^2 +r(s+1)^2=(k+1)s^2+2rs+r.\]
Note that $n^2/(k+1)=(k+1)s^2+2rs+r^2/(k+1)$, so if $k=1$ or 2 this gives $F=\lceil n^2/(k+1) \rceil$. But for $k>2$ the number of filled cells will exceed $\lceil n^2/(k+1) \rceil$ if, as is possible, $r-\lceil r^2/(k+1) \rceil>0$ (see Corollary \ref{cor:main2} for an example).

The construction is illustrated for $(k,n)=(3,16)$ in Figure \ref{fig:3-MOPLS(16)}. This has 64 filled cells and 192 empty cells. It is formed from four copies of a set of three mutually orthogonal Latin squares of order 4 (i.e. four 3-OLS(4)s), placed down the leading diagonal. Here and elsewhere we will write an entry $k$-tuple without commas and brackets, for example $(1,2,3)$ as $123$.

\begin{figure}[ht]
\begin{center}
\setlength\tabcolsep{1.5mm}
\tiny\begin{tabular}{|cccccccccccccccc|}\hline
                \rule{0mm}{3mm}111&222&333&444&-&-&-&-&-&-&-&-&-&-&-&-\\
                432&341&214&123&-&-&-&-&-&-&-&-&-&-&-&-\\
                243&134&421&312&-&-&-&-&-&-&-&-&-&-&-&-\\
                324&413&142&231&-&-&-&-&-&-&-&-&-&-&-&-\\
                -&-&-&-&555&666&777&888&-&-&-&-&-&-&-&-\\
                -&-&-&-&876&785&658&567&-&-&-&-&-&-&-&-\\
                -&-&-&-&687&578&865&756&-&-&-&-&-&-&-&-\\
                -&-&-&-&768&857&586&675&-&-&-&-&-&-&-&-\\
                -&-&-&-&-&-&-&-&999&AAA&BBB&CCC&-&-&-&-\\
                -&-&-&-&-&-&-&-&CBA&BC9&A9C&9AB&-&-&-&-\\
                -&-&-&-&-&-&-&-&ACB&9BC&CA9&B9A&-&-&-&-\\
                -&-&-&-&-&-&-&-&BAC&C9B&9CA&AB9&-&-&-&-\\
                -&-&-&-&-&-&-&-&-&-&-&-&DDD&EEE&FFF&GGG\\
                -&-&-&-&-&-&-&-&-&-&-&-&GFE&FGD&EDG&DEF\\
                -&-&-&-&-&-&-&-&-&-&-&-&EGF&DFG&GED&FDE\\
                -&-&-&-&-&-&-&-&-&-&-&-&FEG&GDF&DGE&EFD\\
                \hline\end{tabular}
                \caption{A 3-MOPLS(16)}\label{fig:3-MOPLS(16)}
                \end{center}
\end{figure}

The $k$-OPLS($n$) produced by Construction 1 is maximal. To see this, let $S_1,S_2,\ldots,S_{k+1}$ denote the $k+1$ entry sets. [In Figure \ref{fig:3-MOPLS(16)},  $S_1=\{1,2,3,4\}$, $S_2=\{5,6,7,8\}$, etc.] If $C$ is an empty cell then, without loss of generality, we may assume that it lies in a row in which all the entries from $S_1$ appear in each of the $k$ coordinate positions, and that it lies in a column in which all the entries from $S_2$ appear in each of the $k$ coordinate positions. Consequently if we wish to place an entry $a_1a_2\ldots a_k$ into cell $C$, none of $a_1,a_2,\ldots, a_k$ can be chosen from $S_1$ or $S_2$. Therefore $k$ choices must be made from the remaining $k-1$ entry sets, and so two or more of these choices must be taken from the same entry set. But every ordered pair from each entry set already appears in the $k$-OPLS($n$), and so it is not possible to place any entry into cell $C$. Hence the array is a $k$-MOPLS($n$).

We conjecture that, for sufficiently large $n$, Construction 1 gives a $k$-MOPLS($n$) with the minimum possible number of filled cells. We also conjecture that any other $k$-MOPLS($n$) with this number of filled cells can be put into the form of Construction 1 by permuting rows and columns of the array. In order to make progress we discuss alternative descriptions for these arrays as graphs and as block designs.

Starting with the simplest case, suppose that $L$ is a Latin square of order $n$. This can be represented as a partition of the complete tripartite graph $K_{n,n,n}$ into edge-disjoint copies of $K_3$. The three vertex sets of the graph represent the rows, columns and entries, and the copies of $K_3$ form the (row, column, entry) triples. Alternatively, $L$ can be represented as a \emph{transversal design}, TD$(3,n)$, which comprises a set of $3n$ points partitioned into three $n$-element subsets, called \emph{groups}, and a set of $n^2$ triples (also called \emph{blocks}) such that each pair of points from different groups appears in precisely one block and no block contains more than one point from each group. The three groups represent the rows, columns and entries, and a block $(a,b,c)$ represents entry $c$ in the cell of $L$ corresponding to the intersection of row $a$ and column $b$.

The other arrays defined above have similar representations. So a \linebreak $k$-OPLS($n$) can be represented as a partition of a subgraph of the complete $(k+2)$-partite graph $K_{(k+2)\times n}$ into copies of $K_{k+2}$. Alternatively, the array can be represented as a partial transversal design, PTD$(k+2,n)$, which comprises a set of $(k+2)n$ points partitioned into $k+2$ ~$n$-element subsets (the groups), and a set of blocks ($(k+2)$-tuples) such that each pair of points from different groups appears in at most one block and no block contains more than one point from each group.

We will switch our viewpoint from arrays to graphs to block designs as the need arises. An advantage of thinking about these arrays in graphical or block design terms is that the vertex sets (or groups) clearly have equal status. So if no special property is assumed for any one of the vertex sets (or groups), then any results for any one of them will, by symmetry, apply to the others.

Suppose that $M$ is a $k$-MOPLS($n$) represented by a graph $G$ as described above. Then the corresponding complementary graph $G'=K_{(k+2)\times n}-G$ is $K_{k+2}$-free. To see this suppose that the vertex sets are $X_i$ for $i=1,2,\ldots,$ $ k+2$, where $X_1$ represents rows and $X_2$ represents columns. If there were a copy of $K_{k+2}$ in $G'$ on the vertex set $\{x_1,x_2,\ldots, x_{k+2}\}$, where $x_i\in X_i$, then cell $(x_1,x_2)$ in the array $M$ could be filled with the entry $k$-tuple $x_3x_4\ldots x_{k+2}$, contradicting the maximality of $M$. The problem of finding the minimum size (number of filled cells) of a $k$-MOPLS($n$) is therefore equivalent to finding the maximum size of the corresponding $K_{k+2}$\,-free complementary graph $G'$.

If there are $F$ filled cells in a $k$-OPLS($n$) then the number of edges between vertex sets $X_i$ and $X_j$ for $1\le i<j\le k+2$ in the corresponding complementary graph $G'$ is $n^2-F$, so that $G'$ has $\binom{k+2}{2}(n^2-F)$ edges. Furthermore, given a vertex $x_i\in X_i$ with degree $d$, the number of edges joining $x_i$ to vertices in $X_j$ ~$(j\ne i)$ is $d/(k+1)$.

Determination of the maximum number of edges in a $K_r$-free graph is a problem that originates in a 1941 paper of Tur\'{a}n's \cite{Turan}, and has been studied extensively in the literature.  In \cite{Bollobas} the authors show that, in particular, the number of edges in a  $K_4$-free $4$-partite graph with parts of size $n$ is at most $5n^2$. Jin \cite{Jin} proved that if $G$ is a $4$-partite graph with parts of size $n$ then $G$ contains a $K_4$ if the minimum vertex degree $\delta(G)$ satisfies $\delta(G)>\lfloor(2+\frac{1}{3})n\rfloor$. Maybe the most important asymptotic contribution to the area is due to Pfender \cite{Pfender} who considers $l$-partite graphs satisfying a certain density condition. For such a graph $G=(X_1\cup X_2\cup ...\cup X_l, E)$, define the density $d(X_i,X_j)$ between the vertex sets $X_i$ and $X_j$ ($i\ne j$) as
\[d(X_i, X_j) = \frac {||G[X_i \cup X_j ]||}{|X_i|\cdot|X_j |},\]
where $||G[X_i \cup X_j ]||$ is the number of edges between $X_i$ and $X_j$, and $|X_h|$ is the number of vertices in $X_h$. Pfender proves that for large enough $l$, and apart from one specified family of graphs, such a graph having $d(X_i, X_j)\ge (r-2)/(r-1)$ for all $i\ne j$ must contain a copy of $K_r$. However, $l$ depends on $r$ and the assumption in the paper \cite{Pfender} is that $l$ is very much larger than $r$. The smallest value of $l$ for which the theorem remains true has not generally been determined.

In the current paper we treat particular instances of the cases $l=r=4$, $l=r=5$, and $l=r\ge 6$. For the complementary graph $G'$ corresponding to a 2-OPLS($n$) with $F$ filled cells, $\displaystyle d(X_i, X_j)=\frac {||G'[X_i \cup X_j ]||}{|X_i|\cdot|X_j |}=\frac{n^2-F}{n^2}$ and we will show that if $F<n^2/3$, so that $d(X_i, X_j)>2/3=(r-2)/(r-1)$, then $G'$ contains a copy of $K_4$ and the 2-OPLS($n$) is not maximal. Similarly for a 3-OPLS($n$), if $F<n^2/4$ so that $d(X_i, X_j)>3/4=(r-2)/(r-1)$, then $G'$ contains a copy of $K_5$ and the 3-OPLS($n$) is not maximal.

\begin{figure}[!ht]\begin{center}
\begin{tikzpicture}[x=0.3mm,y=0.3mm]

\draw (10,10)--(10,110)--(110,110)--(110,10)--(10,10);
\draw [dotted] (10,60)--(110,60);
\draw [dotted] (60,10)--(60,110);

\node at (50,100) {$A$};
\node at (95,100) {$B$};
\node at (50,50) {$C$};
\node at (95,50) {$D$};

\end{tikzpicture}
\caption{A partitioned array. \label{fig:ABCDg}}
\end{center}
\end{figure}
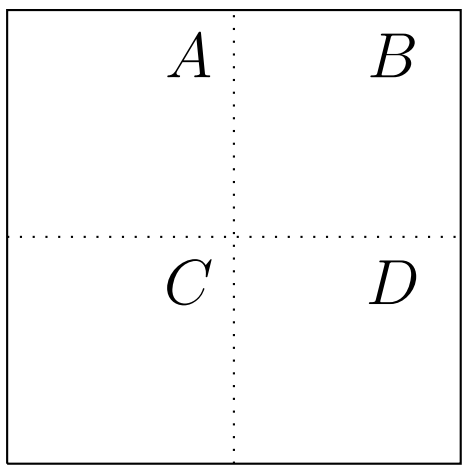

In the proofs we make use of the generic diagram shown in Figure \ref{fig:ABCDg}. This shows a square representing an array that is partitioned into four regions labelled $A,B,C$ and $D$, where $A$ is top left and square, $D$ is bottom right and square, $B$ is top right and $C$ is bottom left. Regions $B$ and $C$ are rectangular, not necessarily square. We will specify the dimensions of these regions in each case as the need arises.

\section{2-MOPLS}\label{sec:2-MOPLS}
For $n\ge 21$, Construction 1 gives a 2-MOPLS($n$) with $\lceil n^2/3 \rceil$ filled cells. To see this, suppose that $n=3s+r$, where $r=0,1 $ or 2, so that $s\ge 7$. For $s\ne 2, 6$ there exists a 2-OLS($s$)(see \cite{HBK}). Hence the construction is always possible for $n\ge 21$. The number of filled cells is $F$ where
\[F=\left\{\begin{array}{ccccc}
3s^2&=&n^2/3 &\mbox{if} &r=0,\\
2s^2+(s+1)^2&=&(n^2 + 2)/3 &\mbox{if} &r=1,\\
s^2+2(s+1)^2&=&(n^2+2)/3 &\mbox{if} &r=2.
\end{array}\right\}=\left\lceil \frac{n^2}{3} \right\rceil. \]
We will prove that these values are best (i.e. minimum) possible. The main step is to
prove Theorem \ref{th:main} below.

Here and elsewhere we use $f(x)$ to denote the \emph{frequency} of a point $x$ in the graphical or block design representation of a $k$-OPLS($n$), that is the number of complete subgraphs ($K_{k+2}$\ s) or blocks ($(k+2)$-tuples) containing $x$. In an array, if $x$ is a row or column then $f(x)$ is the number of entries in that row or column, while if $x$ is an entry then $f(x)$ is the number of times $x$ appears in the array. Note that in the graphical representation of a $k$-OPLS($n$), $f(x)$ is not the degree of the vertex $x$; the degree of $x$ is $(k+1)f(x)$. Of particular interest is the minimum frequency $m$ defined by
\[m=\min\{f(x): x \mbox{~is any vertex}\}.\]

\begin{theorem}\label{th:main}
 Suppose that $M$ is a $2$-MOPLS($n$) with $F$ filled cells. Then $F\ge n^2/3$, with equality if and and only if every frequency (rows, columns, first entries, second entries) is $n/3$.
\end{theorem}

In order to prove this we first establish Lemma \ref{lem:transversals}. To distinguish between rows and columns we will use $r_i$ to denote the row numbered $i$ and $c_i$ to denote the column numbered $i$.

\begin{lemma}\label{lem:transversals}
Suppose that $\Delta$ is a $\delta\times \delta$ array in which some cells may be filled and some may be empty. Let $T$ be a maximum partial transversal of $t$ empty cells in $\Delta$, so that $0\le t\le \delta$. By permuting rows and columns we can assume that $T=\{(r_1,c_1), (r_2,c_2),\ldots, (r_t,c_t)\}$, and let $T'$ be the continuation of this leading diagonal, comprising cells $(r_i,c_i)$ for $t+1\le i\le \delta$. Then every cell $(r_i,c_j)$ with $i,j\in\{t+1, t+2, \ldots,\delta\}$ is filled and $f(r_i)+f(c_j)\ge 2\delta-t$. In particular, for each cell $(r_i,c_i)\in T'$, $f(r_i)+f(c_i)\ge 2\delta-t$.
\end{lemma}
\begin{proof} By saying that $T$ is a maximum partial transversal of empty cells in $\Delta$ we mean that no row or column of $\Delta$ has more than one cell from $T$, and that there is no larger partial transversal of empty cells. We do not exclude the possibilities that $T$ is empty (so that $\Delta$ is full) or that $T$ is a complete transversal of $\delta$ cells.

After appropriate row and column permutations, the array $\Delta$ may be partitioned into four regions labelled $A,B,C,D$ as shown in Figure \ref{fig:ABCDg}. Initially we assume that $t\ne 0,\delta$. Then region $A$ has dimension $t\times t$ and region $D$ has dimension $(\delta-t)\times(\delta-t)$. Region $A$ has the transversal $T$ down its leading diagonal, and region $D$ has the continuation of filled cells $T'$ down its leading diagonal.

The first observation is that region $D$ cannot contain any empty cells since otherwise $T$ would not be a maximum partial transversal of empty cells. So every row and column in the sub-array $D$ has $\delta-t$ filled cells.

Now consider any cell $(r_i,c_j)\in D$, so that $t+1\le i,j\le \delta$, and take $k\in\{1,2,\ldots,t\}$. If both cells $(r_i,c_k)$ and $(r_k,c_j)$ (in regions $C$ and $B$ respectively) were empty then $T$ could be extended by removing cell $(r_k,c_k)$ and adding cells $(r_i,c_k)$ and $(r_k,c_j)$, giving a transversal $\overline{T}$ having $t+1$ empty cells, a contradiction. So at least one of the cells $(r_i,c_k)$ and $(r_k,c_j)$ must be filled for each $k\in\{1,2,\ldots,t\}$. Thus the union of row $r_i$ in region $C$ and column $c_j$ in region $B$ contains at least $t$ filled cells. Since all the cells in region $D$ are filled, we obtain
\[f(r_i)+f(c_j)\ge t+2(\delta-t)=2\delta-t. \]
In particular, this applies when the cell $(r_i,c_j)$ lies in $T'$, i.e. when $i=j$.

In the extreme case $t=0$, region $D$ comprises the whole of the array $\Delta$ and all cells must be filled. At the other extreme $t=\delta$, region $A$ comprises the whole of the array $\Delta$ and there is nothing to prove.
\end{proof}

\begin{proof}[\textit{Proof of Theorem \ref{th:main}}]
Let $m$ denote the minimum frequency in $M$. Then $m\le F/n$. By symmetry we can assume, without loss of generality, that the frequency $m$ is achieved by a first entry $x_1$ and that $x_1$ appears paired with $m$ second entries that we can denote as $y_1,y_2,\ldots,y_m$. If $m>0$, the rows and columns of $M$ can be permuted to partition $M$ as shown in Figure \ref{fig:ABCDg} with entries $x_1y_i$ for $i=1,2,\ldots,m$ down the leading diagonal of region $A$, so that $A$ has dimension $m\times m$ and region $D$ has dimension $(n-m)\times (n-m)$. There is no entry $x_1$ in any column through region $B\cup D$ and in any row through region $C\cup D$. If $m=0$ then regions $A,B,C$ do not exist, and region $D$ comprises the whole of $M$.

If $I=(r_i,c_i)$ is one of the $m$ cells on the leading diagonal in region $A$, then $f(r_i)\ge m$ and $f(c_i)\ge m$, so $f(r_i)+f(c_i)\ge 2m$.

Region $D$ will contain some filled cells and some empty cells. We may apply Lemma \ref{lem:transversals} to region $D$ by taking  $\delta=n-m$.

If $I=(r_i,c_i)$ is one of the $t$ empty cells in the maximum partial transversal $T$ of empty cells in $D$, then the union of row $r_i$ and column $c_i$ must contain at least $n-m$ entries, otherwise there would be some second entry $y_j$ with $j\ne 1,2,\ldots,m$ that could be paired with $x_1$ and placed in $I$, contradicting the maximality of $M$. Hence $f(r_i)+f(c_i)\ge n-m$ for each $i\in\{m+1,$ $m+2,\ldots,m+t\}$.

If $I=(r_i,c_i)$ is one of the $n-m-t$ filled cells in region $D$ lying in the continuation partial transversal $T'$ then, by Lemma \ref{lem:transversals} $f(r_i)+f(c_i)\ge 2(n-m)-t$.

Summing over all the cells on the leading diagonal of $M$ gives
\[\sum_{i=1}^n \left(f(r_i)+f(c_i)\right)\ge 2m^2+t(n-m)+(n-m-t)(2(n-m)-t).\]

Since $\displaystyle \sum_{i=1}^n f(r_i)=F=\sum_{i=1}^n f(c_i),$ we obtain
\begin{align*}2F &\ge2m^2+(n-m)^2+(n-m-t)^2=(n-m-t)^2+\frac{(n-3m)^2}{3}+\frac{2n^2}{3}.
\end{align*}
This gives
\begin{align}F&\ge \left\lceil\frac{(n-m-t)^2}{2}+\frac{(n-3m)^2}{6}+\frac{n^2}{3}\right\rceil.\label{ineq:main}\end{align}

It follows that $F\ge n^2/3$. From the inequality (\ref{ineq:main}) it is clear that the equality $F=n^2/3$ is only possible if $t=n-m$, and $m=n/3$, which of course implies that $n\equiv 0$ (mod 3). And if $F=n^2/3$ and $m=n/3$ then the frequency of every row (also every column, first entry, second entry) is $n/3$. Conversely, if the frequency of every row is $n/3$, then $\displaystyle F=\sum_{i=1}^n f(r_i)=n^2/3$.
\end{proof}

\noindent\textbf{Further consequences}\\
We can extract more information from inequality (\ref{ineq:main}). If $n=3s+1$ it gives
\[F\ge \left\lceil\frac{9s^2+6s+1}{3}\right\rceil= \frac{9s^2+6s+3}{3}=\frac{n^2+2}{3}.\]
If $n=3s+2$ the inequality gives
\[F\ge \left\lceil\frac{9s^2+12s+4}{3}\right\rceil= \frac{9s^2+12s+6}{3}=\frac{n^2+2}{3}.\]
In each case equality requires that $m=s=\lfloor n/3\rfloor$. These observations, taken with Construction 1, enable us to state the following corollary.

\begin{cor}\label{cor:main}
For $n\ge 21$, the minimum number of filled cells $F$ achievable in a $2$-MOPLS($n$) is given by
\[F=\left\{\begin{array}{cl}
n^2/3 &\mbox{if~~}  n\equiv 0 \pmod{3}\\
(n^2 + 2)/3 &\mbox{if~~} n\equiv 1 \mbox{~or~} 2 \pmod{3}
\end{array}\right\} = \left\lceil\frac{n^2}{3}\right\rceil.\]
Moreover, the minimum frequency in each of these cases is $m=\lfloor n/3 \rfloor$.
\end{cor}

\section{3-MOPLS}\label{sec:3-MOPLS}
For $n\ge 44$, Construction 1 gives a 3-MOPLS($n$) with approximately $n^2/4$ filled cells. To see this, suppose that $n=4s+r$, where $r=0,1,2$ or 3, so that $s\ge 11$. For $s\ne 2,3,6,10$ there exists a 3-OLS($s$) (see \cite{HBK}). Hence the construction is always possible for $n\ge 44$. The number of filled cells is $F$ where
\[F=\left\{\begin{array}{ccccc}
4s^2&=&n^2/4 &\mbox{if} &r=0,\\
3s^2+(s+1)^2&=&(n^2 + 3)/4 &\mbox{if} &r=1,\\
2s^2+2(s+1)^2&=&n^2/4 +1&\mbox{if} &r=2,\\
s^2+3(s+1)^2&=&(n^2+3)/4&\mbox{if} &r=3.
\end{array}\right. \]
We will prove that these values are best (i.e. minimum) possible. We can regard a 3-MOPLS($n$) either as an array with triples in occupied cells, or as a partial transversal design PTD($5,n$) having five groups, $X_i$ for $i=1,2,3,4,5$ with blocks $x_1x_2x_3x_4x_5$ ($x_i\in X_i$) corresponding to the entry triple $x_3x_4x_5$ in cell $(x_1,x_2)$ of $M$.

The first step is to establish an upper bound on the combined frequencies of a non-occurring pair of points in a 3-MOPLS($n$) that has $F$ filled cells. We say that points $u$ and $v$ form a \emph{non-occurring pair of points} if $u$ and $v$ are in different groups and there is no block of the design that contains both points. The result actually applies to $k$-MOPLS($n$) for any $k\ge 2$.

\begin{lemma}\label{lemma:ab}
 Suppose that $M$ is a $k$-MOPLS($n$) with $F<n^2$ filled cells. Then for any pair of groups $(X_i,X_j)$ ~($i\ne j$) there exists a non-occurring pair of points $ab$ with $a\in X_i$, $b\in X_j$, and $f(a)+f(b)\le 2F/n$. Unless every frequency in that pair of groups is $F/n$, there is a pair of non-occurring points $ab$ in that pair of groups with $f(a)+f(b)<2F/n$.
\end{lemma}
\begin{proof}
 We calculate the average of $f(r)+f(c)$ over all non-occurring pairs $rc$, where $r$ is a row and $c$ is a column (i.e. cell $(r,c)$ is empty). Denote the set of all these pairs by $\mathcal{E}$, so that $|\mathcal{E}|=n^2-F$. In the summation $\displaystyle\sum_{rc\in\mathcal{E}}\left(f(r)+f(c)\right)$, ~$f(r)$ will be counted $(n-f(r))$ times.  Hence
 \[\sum_{rc\in\mathcal{E}}f(r)=\sum_{r=1}^n \big(f(r)(n-f(r))\big)= nF-\sum_{r=1}^n (f(r))^2. \]
 But
 \[\sum_{r=1}^n (f(r))^2 \ge \frac1n \left(\sum_{r=1}^n f(r)\right)^2 =\frac{F^2}{n},\]
 with equality only if $f(r)=F/n$ for $r=1,2,\ldots,n$.
 Applying the same reasoning to columns and combining the inequalities gives
\[\sum_{rc\in\mathcal{E}}\left(f(r)+f(c)\right)\le2\left(nF-\frac{F^2}{n}\right) =\frac{2F}{n}\left(n^2-F\right)=\frac{2F}{n}|\mathcal{E}|.\]
It follows that there must be a non-occurring pair $rc$ with $f(r)+f(c)\le 2F/n$, and if any row or column has frequency different from $F/n$ then we have a non-occurring pair $rc$ with $f(r)+f(c)<2F/n$.

Although the argument above singles out rows and columns for special treatment, symmetry between groups of the design ensures that the result applies to every pair of groups. \end{proof}

We can now prove the main result of this Section.
\begin{theorem}\label{th:3-MOPLS}
 Suppose that $M$ is a $3$-MOPLS($n$) with $F$ filled cells. Then $F\ge n^2/4$, with equality if and only if every frequency is $n/4$.
\end{theorem}
\begin{proof} Suppose that $M$ is a 3-MOPLS($n$) with $F$ filled cells. We may assume that $F<n^2/2$. Let $\mu$ denote the minimum value of $f(u)+f(v)$ for non-occurring pairs $uv$. By Lemma \ref{lemma:ab}, $\mu< n$. Suppose that $a$ and $b$ achieve this bound, that is $f(a)+f(b)=\mu$. Without loss of generality, we may assume that $a$ is a first entry and $b$ is a second entry in the array $M$.

If $\mu=0$ then neither $a$ nor $b$ occur in $M$ and so for any third entry $c$ the entry triple $abc$ could be placed in any empty cell $(i,j)$ unless $c$ is present in row $i$ or column $j$. The maximality of $M$ therefore implies that $M$ forms an MPLS($n$) on third entries, and so by \cite{HorakRosa} $M$ has at least $n^2/2$ filled cells, a contradiction. So $\mu\ne 0$.

We may also permute the rows and columns of $M$ so that all occurrences of $a$ and $b$ are in the subsquare formed by rows and columns $1,2,\ldots,\mu$. It is possible that some rows may contain both $a$ and $b$ in different cells, in which case the total number of rows containing either $a$ or $b$ will be less than $\mu$. A similar observation applies to columns.

The rows and columns of $M$ can be rearranged as shown in Figure \ref{fig:3-MOPLS(n)}, where $M$ is partitioned into four main regions labelled $A,B,C,D$.  Region $A$ has dimension $\mu\times\mu$ and $D$ has dimension $(n-\mu)\times(n-\mu)$. Since $\mu<n$, region $D$ must exist.

\begin{figure}[!ht]\begin{center}
\begin{tikzpicture}[x=1mm,y=1mm]

\draw (10,10)--(10,110)--(110,110)--(110,10)--(10,10);
\draw [dotted] (10,70)--(110,70);
\draw [dotted] (50,10)--(50,110);

\draw [dotted] (10,90)--(28,90)--(28,110);
\draw [dotted] (50,44)--(74,44)--(74,70);

\node at (40,105) {$A$};
\node at (100,105) {$B$};
\node at (40,65) {$C$};
\node at (100,65) {$D$};

\node at (12,107) {$\emptyset$};
\node at (19,101) {$\ddots$};
\node at (26,93) {$\emptyset$};
\node at (30,90) {$\star$};
\node at (37,83) {$\ddots$};
\node at (47,72) {$\star$};

\node at (53,67) {$\emptyset$};
\node at (62,57) {$\ddots$};
\node at (72,47) {$\emptyset$};
\node at (76,42) {$\star$};
\node at (88,31) {$\ddots$};
\node at (108,12) {$\star$};

\node at (80,120) {No $a$ or $b$ in};
\node at (80,115) {these columns};
\draw [Latex-] (51,115)--(65,115);
\draw [-Latex] (95,115)--(109,115);
\node at (-3,42.5) {No $a$ or $b$ in};
\node at (-3,37.5) {these rows};
\draw [Latex-] (-3,11)--(-3,33.5);
\draw [-Latex] (-3,46.5)--(-3,69);

\draw [Latex-Latex] (10,8)--(49,8);
\node at (31,3) {$\mu$};
\draw [Latex-Latex] (51,8)--(110,8);
\node at (80,3) {$n-\mu$};

\draw [Latex-Latex] (52,71)--(52,109);
\node at (57,92) {$\mu$};
\draw [Latex-Latex] (48,11)--(48,69);
\node at (40,45) {$n-\mu$};

\draw [Latex-Latex] (7,90)--(7,109);
\node at (4,100) {$p$};

\draw [Latex-Latex] (76,44)--(76,69);
\node at (79,57) {$q$};

\end{tikzpicture}
\caption{A partition of $M$.\label{fig:3-MOPLS(n)}}
\end{center}\vspace{-5mm}
\end{figure}
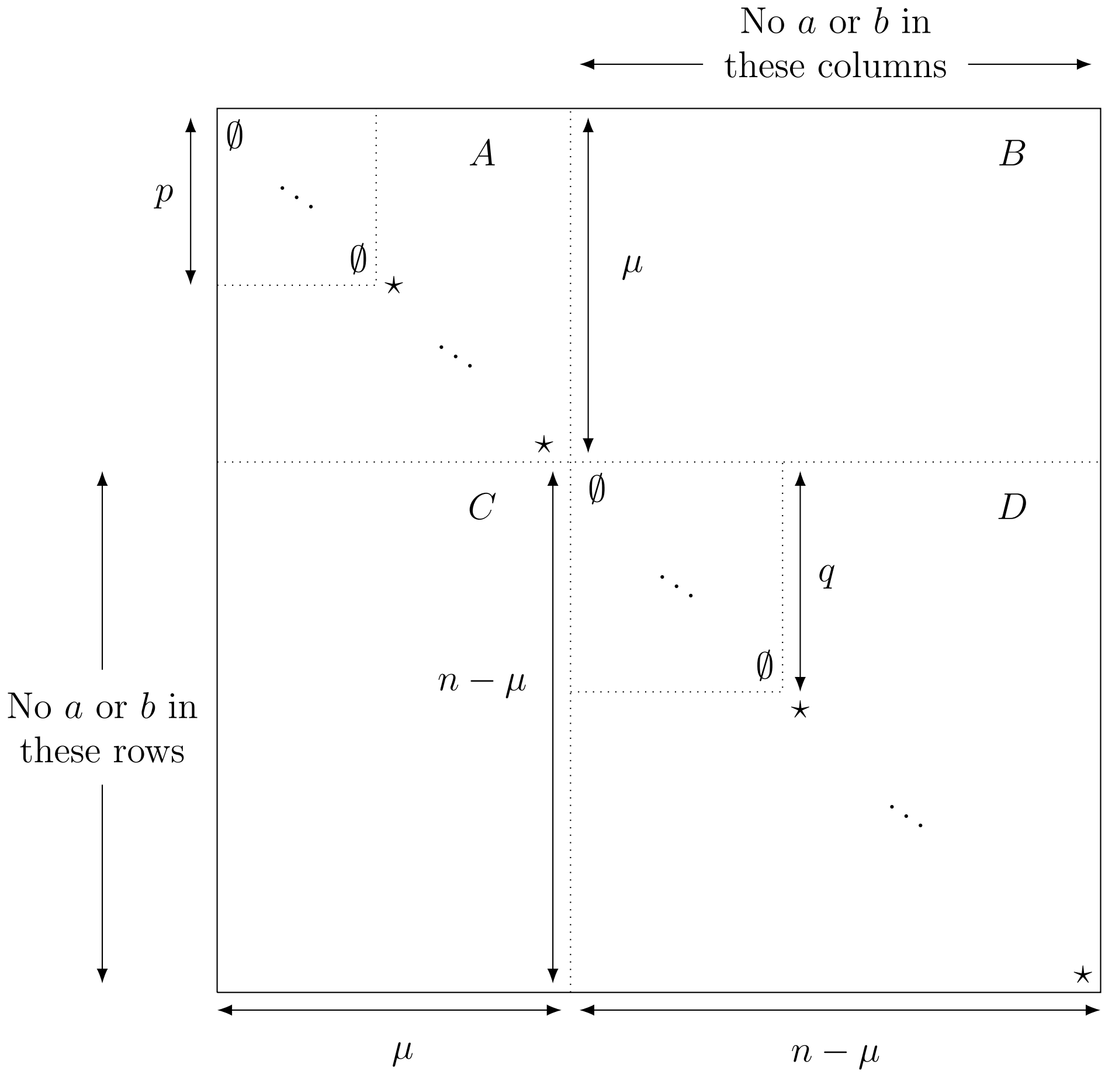

By permuting rows and columns numbered from 1 to $\mu$, a maximum partial transversal of $p\ge 0$ empty cells in region $A$ can be located on the leading diagonal of $A$. These empty cells are indicated in Figure \ref{fig:3-MOPLS(n)} by the symbol $\emptyset$. Any remaining cells on the leading diagonal of region $A$ must be occupied by some triple, and these occupied cells are denoted by the symbol $\star$.

Region $D$ may be treated in a similar way by permuting rows and columns numbered from $\mu+1$ to $n$, so that it has a maximum partial transversal of $q\ge0$ empty cells on its leading diagonal.

With the array $M$ arranged as described, the proof proceeds by estimating the sum of row frequency and column frequency for each cell on the leading diagonal of $M$. To distinguish between rows and columns we again use $r_i$ to denote the row numbered $i$ and $c_i$ to denote the column numbered $i$.

If the cell $(r_i,c_i)$ lies in the first $p$ cells on the leading diagonal of region $A$, then $r_ic_i$ is a non-occurring pair. Hence $f(r_i)+f(c_i)\ge \mu$.

If the cell $(r_i,c_i)$ is one of the remaining cells on the leading diagonal of region $A$, then we may apply Lemma \ref{lem:transversals} on $A$ to give $f(r_i)+f(c_i)\ge 2\mu-p$.

If the cell $(r_i,c_i)$ lies in the first $q$ cells on the leading diagonal of region $D$ then the union of row $r_i$ and column $c_i$ must contain at least $n-f(a)-f(b)=n-\mu$ filled cells. This is because there are exactly $f(a)$ different third entries appearing with $a$ in a block of $M$, and likewise $f(b)$ third entries appearing with $b$. So there are at least $n-f(a)-f(b) ~(>0)$ third entries that do not occur with either $a$ or $b$. If any one of these, say $z$, is not in the union of row $r_i$ and column $c_i$ then the triple $abz$ may be inserted into cell $(r_i,c_i)$, contradicting the maximality of $M$. It follows that for these cells $f(r_i)+f(c_i)\ge n-\mu$.

If the cell $(r_i,c_i)$ is one of the remaining cells on the leading diagonal of region $D$, then we may apply Lemma \ref{lem:transversals} on $D$ to give $f(r_i)+f(c_i)\ge 2(n-\mu)-q$.

Summing over all cells on the leading diagonal of $M$ gives
\begin{align*}2F&=\sum_{i=1}^n\big(f(r_i)+f(c_i)\big)\\
&\ge p\mu+(\mu-p)(2\mu-p)+q(n-\mu)+(n-\mu-q)(2(n-\mu)-q).
\end{align*}
Hence
\begin{align*}
 2F&\ge\mu^2+(\mu-p)^2+(n-\mu)^2+(n-\mu)^2-2q(n-\mu)+q^2\\
 &=\mu^2+(\mu-p)^2+(n-\mu)^2+(n-\mu-q)^2\\
 &=(\mu-p)^2+(n-\mu-q)^2+\frac{(n-2\mu)^2}{2}+\frac{n^2}{2}.
 \end{align*}
This gives
\begin{align}F&\ge \left\lceil \frac{(\mu-p)^2}{2} +\frac{(n-\mu-q)^2}{2} + \frac{(n-2\mu)^2}{4}+ \frac{n^2}{4}\right\rceil \label{ineq:main2}.\end{align}
It follows that $F\ge n^2/4$.

Now suppose that $F=n^2/4$ and that the minimum frequency $m<n/4=F/n$. By Lemma \ref{lemma:ab} there exist points $u,v$ such that $uv$ forms a non-occurring pair and $f(u)+f(v)<2F/n$. Hence $\mu<n/2$ and inequality (\ref{ineq:main2}) gives a contradiction. Thus the only possibility for $F=n^2/4$  is that every point has frequency $n/4$ and $\mu=p=q=n/2$. Consequently $n\equiv 0$ (mod 4) and $M$ has a full transversal of empty cells.
\end{proof}

\noindent\textbf{Further consequences}\\
We can extract more information from inequality (\ref{ineq:main2}). If $n=4s+1$ it gives
\[F\ge \left\lceil\frac{16s^2+8s+1}{4}\right\rceil= \frac{16s^2+8s+4}{4}=\frac{n^2+3}{4}.\]
If $n=4s+2$ then the inequality gives $F\ge n^2/4$, but equality is only possible if $m=n/4$, which is not an integer, and so $F\ge n^2/4 +1$.\\
If $n=4s+3$ the inequality gives
\[F\ge \left\lceil\frac{16s^2+24s+9}{4}\right\rceil= \frac{16s^2+24s+12}{4}=\frac{n^2+3}{4}.\]
These observations, taken with Construction 1, enable us to state the following corollary.

\begin{cor}\label{cor:main2}
For $n\ge 44$, where $n=4s+r$ and $0\le r\le 3$, the minimum number of filled cells $F$ achievable in a $3$-MOPLS($n$) is given by
\[F=\left\{\begin{array}{ccccc}
4s^2&=&n^2/4 &\mbox{if} &r=0,\\
3s^2+(s+1)^2&=&(n^2 + 3)/4 &\mbox{if} &r=1,\\
2s^2+2(s+1)^2&=&n^2/4 +1&\mbox{if} &r=2,\\
s^2+3(s+1)^2&=&(n^2+3)/4&\mbox{if} &r=3.
\end{array}\right. \]
\end{cor}

It is a little more tricky than it was in the case of 2-MOPLS($n$) to identify the minimum frequency in each of these cases. But the following theorem does the job.

\begin{theorem}\label{th:minfreq}
 Suppose that $M$ is a $3$-MOPLS($n$) ($n\ge 44$) with the minimum number of filled cells. Then the minimum frequency $m=\lfloor n/4\rfloor$.
\end{theorem}
\begin{proof}  Theorem \ref{th:3-MOPLS} deals with the case when $n\equiv 0$ (mod 4), so assume that $n=4s+r$ where $r=1,2$ or 3. Corollary \ref{cor:main2} gives the number of filled cells in $M$.

The first step is to prove that if $m\le\lfloor n/4\rfloor-1=s-1$ then there would be a non-occurring pair $ab$ with $f(a)+f(b)\le (n-3)/2$. Denote the five groups of $M$ as $X_i$ for $i=1,2,\ldots,5$. Without loss of generality, we may assume that there is a point $a\in X_1$ with $f(a)=m\le\lfloor n/4\rfloor-1$. There will be $n-m$ points in $X_2$ that do not lie in a block with $a$. Denote this set of points as $X_2'$. Now suppose that each of these points has frequency in excess of $(n-3)/2-m$, i.e. if $x\in X_2'$ then $f(x)>(n-3)/2-m$. So if $r=1$ or 3, $f(x)\ge (n-1)/2-m$, while if $r=2$, $f(x)\ge (n-2)/2-m$.

Counting frequencies of all the points in $X_2$ gives
\[F\ge m^2 +(n-m)\left(\frac{n-\alpha}{2}-m\right),\]
where $\alpha=1$ if $r=1$ or 3, while $\alpha=2$ if $r=2$. Arithmetic then gives
\[F\ge \frac{n^2}{4} +\frac{(n-4m)(n-2m)}{4}-\frac{\alpha}{2}(n-m).\]

If $r=1$ and $m\le s-1$, then $n-4m\ge 5$ and so
\[F\ge \frac{n^2}{4}+\frac{5(n-2m)}{4}-\frac{(n-m)}{2}=\frac{n^2}{4}+\frac{3n-8m}{4}\ge \frac{n^2}{4}+\frac{n+10}{4}>\frac{n^2+3}{4},\]
a contradiction.

If $r=2$ and $m\le s-1$, then $n-4m\ge 6$ and so
\[F\ge \frac{n^2}{4}+\frac{6(n-2m)}{4}-\frac{2(n-m)}{2}=\frac{n^2}{4}+\frac{2n-8m}{4}\ge \frac{n^2}{4}+3>\frac{n^2+4}{4},\]
a contradiction.

If $r=3$ and $m\le s-1$, then $n-4m\ge 7$ and so
\[F\ge \frac{n^2}{4}+\frac{7(n-2m)}{4}-\frac{(n-m)}{2}=\frac{n^2}{4}+\frac{5n-12m}{4}\ge \frac{n^2}{4}+\frac{2n+21}{4}>\frac{n^2+3}{4},\]
a contradiction.

In each case the contradiction implies that there must be a point $b\in X_2'$ such that $f(b)\le (n-3)/2-m$ and so $f(a)+f(b)\le (n-3)/2$. Hence $\mu\le (n-3)/2$. Inequality (\ref{ineq:main2}) from the proof of Theorem \ref{th:3-MOPLS} gives
\begin{align*}F&\ge \left\lceil \frac{(n-2\mu)^2}{4}+ \frac{n^2}{4}\right\rceil
 \ge \left\lceil \frac94+ \frac{n^2}{4}\right\rceil>\frac{n^2+4}{4},
\end{align*}
which is a further contradiction. Hence $m\ge\lfloor n/4\rfloor$, and since $m$ cannot exceed the mean frequency, namely $F/n$, it follows that $m=\lfloor n/4\rfloor$.
\end{proof}

In the next section we prove some results about the structure of a \linebreak $k$-MOPLS($n$) in certain cases. These will confirm that, for sufficiently large $n$, any minimal 2-MOPLS($n$) and any minimal 3-MOPLS($n$) essentially have the form produced by Construction 1.

\section{\textit{k}-MOPLS}\label{sec:k-MOPLS}
Theorem \ref{th:main} resolves the main conjecture concerning 2-MOPLS made in \cite{MOPLS1}, while Theorem \ref{th:3-MOPLS} deals with the extension of this to the case of 3-MOPLS. We are still unable to prove that, for $k\ge 4$, Construction 1 always gives the minimum number of filled cells in a $k$-MOPLS($n$) when $n$ is sufficiently large. However, we can prove some results about the structure of a $k$-MOPLS($n$) in certain cases. We start with the case $k=1$; that is the equivalent result for MPLS($n$)s. Horak and Rosa \cite{HorakRosa} proved that the minimum number of filled cells in an MPLS($n$) is $\lceil n^2/2 \rceil$.

\begin{lemma} \label{lem:HR} Suppose that $n>1$ and that $M$ is an MPLS($n$) having exactly $F=\lceil n^2/2 \rceil$ filled cells. Then $M$ can be partitioned into two Latin squares, on different point sets and occupying different rows and columns, with orders that sum to $n$ and are either $\lfloor n/2 \rfloor$ or $\lfloor n/2 \rfloor+1$.
\end{lemma}
\begin{proof}
If we put $s=\lfloor n/2 \rfloor$, then $n=2s$ or $2s+1$, so that $M$ has $F=2s^2$ or $2s^2+2s+1$ filled cells, respectively. Let $m$ be the minimum frequency of $M$.  Then $m$ is at most the mean frequency, namely $F/n$, which gives $m\le s$ in both cases. By the symmetry between rows, columns and entries, we may assume that there is some entry, say $x_1$, that appears in $M$ exactly $m$ times and, by row and column permutations, $M$ may then be rearranged into the form shown in Figure \ref{fig:ABCDg}, where $x_1$ appears only in region $A$ and fills its main diagonal. Then region $A$ has dimension $m\times m$, and consequently region $D$ has dimension $(n-m)\times(n-m)$.

Region $D$ must be full, otherwise $x_1$ could be inserted into an empty cell contradicting the maximality of $M$. Since $m$ is the minimum frequency we then see that $M$ has $m$ rows each with at least $m$ entries and $n-m$ rows each with at least $n-m$ entries, giving $F\ge m^2+(n-m)^2=2m^2-2mn +n^2$. If $n=2s$ this gives $2s^2\ge 2m^2-4ms+4s^2=2s^2+2(s-m)^2$, which is only possible if $m=s$. Similarly, if $n=2s+1$, we get $2s^2+2s+1\ge 2m^2-4ms-2m+4s^2+4s+1$, and this gives $2s^2\ge 2s^2+2(s-m)^2+2(s-m)$, which again is only possible if $m=s$. In either case $F=m^2+(n-m)^2$, and it then follows that $M$ has $m$ rows with exactly $m$ entries and  $n-m$ rows with exactly $n-m$ entries. By symmetry between rows and columns, the same must be true for columns. Consequently regions $B$ and $C$ must be empty and regions $A$ and $D$ must be full.

Let $S$ denote the set of $n-m$ entries occurring in column $m+1$ of $M$, all of which occur in region $D$. We may take $S=\{x_{m+1},x_{m+2},\ldots, x_n\}$. Since $M$ is maximal, each row $i$ of $M$ for $1\le i\le m$ must contain all the entries in $S'=\{x_1,x_2,\ldots,x_m\}$, otherwise there exists some $k\le m$ such that the entry $x_k$ can be placed in cell $(i,m+1)$. Consequently each element of $S'$ must occur in each of the rows $1,2,\ldots,m$ of $M$. Since $M$ is a partial Latin square, it follows that $A$ is a Latin square of order $m$ with entries from $S'$. Next consider the set of entries occurring in row $1$ of $M$. By employing the same argument again it can be seen that each column $m+i$ of $M$ for $1\le i \le n-m$ must contain all the entries in $S$. Hence $D$ is a Latin square of order $n-m$ with entries from $S$.
\end{proof}

Using Lemma \ref{lem:HR} the main result of this Section can be established.

\begin{theorem}\label{th:partition}
Suppose that $M$ is a $k$-MOPLS($n$) with $n=(k+1)s+r$, where $0\le r\le k$, having exactly $F= (k+1-r)s^2+r(s+1)^2$ filled cells, and with all frequencies either $s~(=\lfloor n/(k+1)\rfloor)$ or $s+1$. Then $M$ can be partitioned into $k+1-r$ $k$-OLS($s$)s and $r$\,$k$-OLS($s+1$)s on different point sets and occupying different rows and columns.
\end{theorem}
\begin{proof}
The conditions on $F$ and the frequencies ensure that $M$ has the same parameters (number of filled cells and frequencies) as the $k$-MOPLS($n$) produced by Construction 1. The theorem says that, essentially, $M$ must have the form produced by the construction. The proof is by induction on $k$. For the purpose of the proof, an MPLS($n$) can be regarded as a 1-MOPLS($n$). Lemma \ref{lem:HR} establishes the result for $k=1$.

So, as an inductive assumption, suppose that $k\ge 2$ and that any ($k-1$)-MOPLS($n$) with $n=ks+r$, where $0\le r\le k-1$, having exactly $F= (k-r)s^2+r(s+1)^2$ filled cells, and with all frequencies either $s$ or $s+1$, can be partitioned as described.

Viewing $M$ as a block design, consider any group $X_i$. Let $n_0$ denote the number of points $x_i\in X_i$ with $f(x_i)=s$, and let $n_1$ denote the number of points $x_i\in X_i$ with $f(x_i)=s+1$.
Then
\[n_0+n_1=n=(k+1)s+r, \mbox{~and~~}sn_0 +(s+1)n_1=F=(k+1-r)s^2+r(s+1)^2 .\]
Hence
\[n_0 =(k+1-r)s, \mbox{~and~~}n_1 =r(s+1).\]
These equalities hold for every group $X_i$, $i=1,2,\ldots,k+2$.

Now let $X_{k+1}$ represent rows and $X_{k+2}$ represent columns. We will select a $k$-tuple of entries, containing only non-occurring pairs. First choose $a_1\in X_1$ with $f(a_1)=s$. Then there are $n-s$ choices for $a_2\in X_2$ such that $a_1a_2$ is a non-occurring pair. At least $n_0-s$ of these choices have frequency $s$. So take $a_2$ with $f(a_2)=s$.

Now there are $n-2s$ choices for $a_3\in X_3$ such that $a_1a_3$ and $a_2a_3$ are non-occurring pairs. At least $n_0-2s$ of these choices have frequency $s$. So take $a_3$ with $f(a_3)=s$.

We can proceed in this fashion to take $a_i\in X_i$ with $f(a_i)=s$ until $i$ is so large that $n_0-(i-1)s\le 0$, equivalently $i\ge k+2-r$. So for $i=1,2,\ldots,k+1-r$ ($r\ge 1$) it is possible to choose $a_i$ such that $f(a_i)=s$ and $a_1a_2\ldots a_i$ has only non-occurring pairs. For larger values of $i$ (i.e. $i=k+2-r,k+3-r,\ldots,k$, where $r\ge 2$) we only have $f(a_i)\le s+1$, and there are at least $n-(i-1)(s+1)+(k+1-r)=(k+2-i)(s+1)>0$ choices for $a_i$ such that $a_1a_2\ldots a_i$ has only non-occurring pairs. It follows that for $r\ge1$
\[\sum_{i=1}^k f(a_i)\le (k+1-r)s+(r-1)(s+1)=n-(s+1).\]
If $r=0$ then each $a_i$ is chosen with frequency $s$ for $i=1,2,\ldots, k$, and then $\sum_{i=1}^k f(a_i)=ks=n-s$. This case when $r=0$ will be dealt with separately.

\medskip

Now assume that $r\ge 1$. Note that if any $a_i$ for $i> k+1-r$ is chosen with $f(a_i)=s$  then $\sum_{i=1}^k f(a_i)<n-(s+1)$. So equality can only happen if every $a_i\in X_i$ for $i>k+1-r$ has $f(a_i)=s+1$. In general, the inequality $\sum_{i=1}^k f(a_i)\le n-(s+1)$ implies that $M$ may be partitioned as shown in Figure \ref{fig:ABCDg} with region $A$ having dimension $\sum_{i=1}^k f(a_i) \times \sum_{i=1}^k f(a_i)$ and region $D$ having dimension $d\times d$, where $d\ge s+1$. All occurrences of $a_1,a_2,\ldots,a_k$ are in region $A$, so region $D$ must be full, otherwise $a_1a_2\ldots a_k$ could be inserted into an empty cell. This implies $d\le s+1$, so that $d=s+1$. But then regions $B$ and $C$ must be empty since rows and columns have frequency at most $s+1$. It also follows that $\sum_{i=1}^k f(a_i)=n-(s+1)$, and so $f(a_i)=s$ for $i=1,2,\ldots, k+1-r$ and $f(a_i)=s+1$ for $i=k+2-r,k+3-r,\ldots,k$.

If some row through region $A$ contained more than one of $a_1,a_2,\ldots,a_k$ in different cells, then at least one row through region $A$ would not contain any of these entries. This in turn would imply that region $A$ could be reduced by at least one row, with a corresponding increase of at least one row in region $D$. But then each column of region $D$ would have at least $s+2$ entries, which contradicts the assumption that each column has frequency at most $s+1$. So each row through region $A$ contains precisely one of $a_1,a_2,\ldots,a_k$, with a similar result for columns.  The rows and columns through region $A$ may therefore be permuted so that the $a_i$\,s occur in order on the leading diagonal of $A$ and, for each $i$, the $s$ or $s+1$ occurrences of $a_i$ are enclosed in a square subregion denoted by $A_i$ that has dimension either $s\times s$ or $(s+1)\times (s+1)$. This is shown in Figure \ref{fig:diag2}.

\begin{figure}[ht!]
\begin{center}
\begin{tikzpicture}[x=0.5mm,y=0.5mm]

\draw (10,10)--(10,110)--(110,110)--(110,10)--(10,10);
\draw [dotted] (10,30)--(110,30);
\draw [dotted] (90,10)--(90,110);

\node at (70,100) {$A$};
\node at (100,100) {$B$};
\node at (70,20) {$C$};
\node at (100,20) {$D$};

\draw [dotted] (10,90)--(30,90)--(30,110);
\node at (20,100){$A_1$};
\draw [dotted] (30,90)--(50,90)--(50,70)--(30,70)--(30,90);
\node at (40,80){$A_2$};
\node at (60,60){$\ddots$};
\draw [dotted] (70,30)--(70,50)--(90,50);
\node at (80,40){$A_k$};

\draw [Latex-Latex] (10,5)--(89,5);
\node at (50,-5) {$n-(s+1)$};
\draw [Latex-Latex] (91,5)--(110,5);
\node at (100,-5) {$s+1$};

\draw [Latex-Latex] (5,31)--(5,109);
\node at (-15,70) {$n-(s+1)$};
\draw [Latex-Latex] (5,11)--(5,29);
\node at (-12,20) {$s+1$};

\node at (160,90) {\small Each $A_i$ has the form};
\draw (140,40)--(140,80)--(180,80)--(180,40)--(140,40);
\node at (145,75){$a_i$};
\node at (160,60){$\ddots$};
\node at (175,45){$a_i$};
\node at (175,30) {\small Dimension $s\times s$ if $i\le k+1-r$,};
\node at (175,20){\small$(s+1)\times(s+1)$ otherwise.};

\end{tikzpicture}
\caption{A finer partition of $M$.\label{fig:diag2}}
\end{center}
\end{figure}
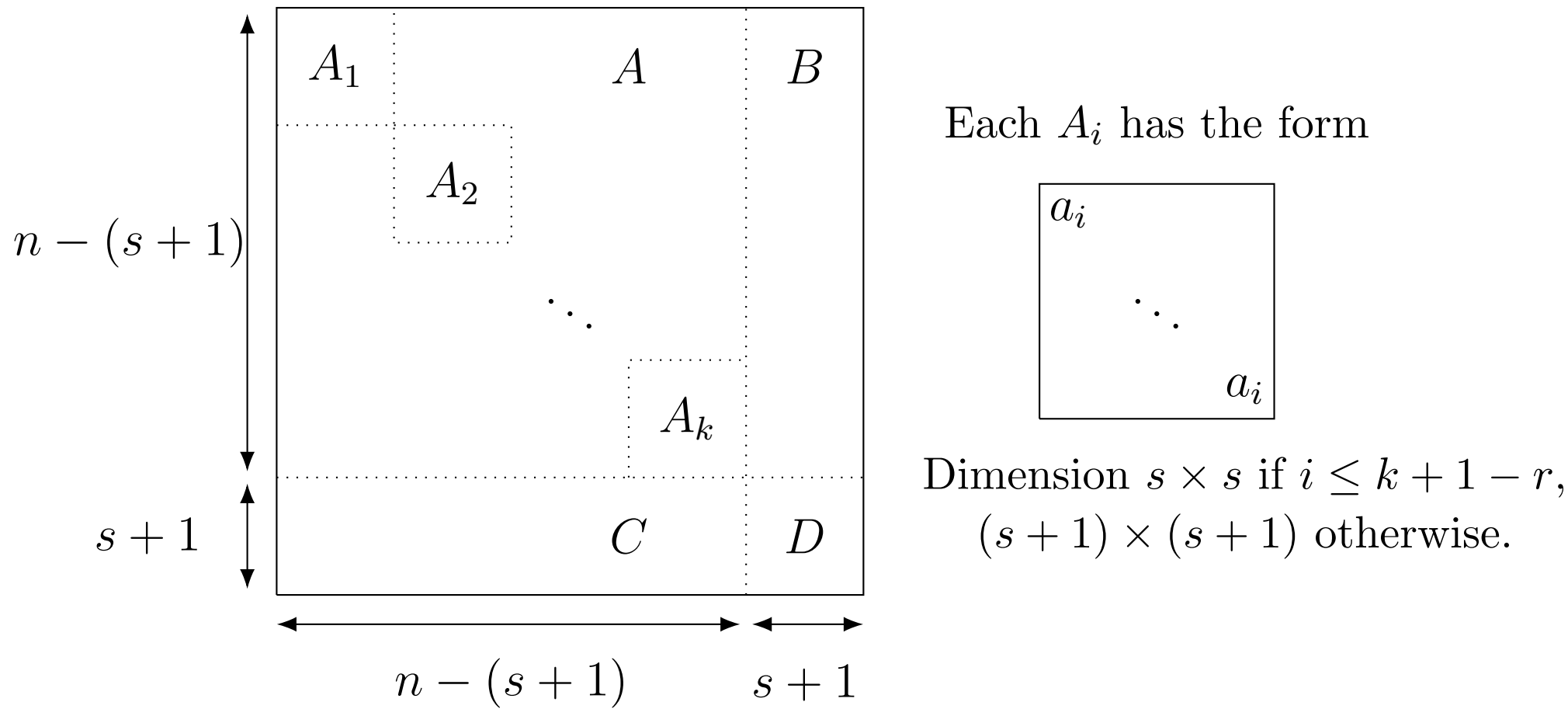

We will prove that $D$ forms a Latin square on the $i^{\mathrm{th}}$ coordinate for each $i=1,2,\ldots,k$. There are two cases to consider: (a) when $A_i$ has dimension $s\times s$, and (b) when $A_i$ has dimension $(s+1)\times(s+1)$.

\medskip

\noindent\textbf{Case (a): $A_i$ has dimension $s\times s$.} Consider the $i^{\mathrm{th}}$ coordinate entries that do not occur in a block with any $a_j$ for $j\ne i$. There must be at least $2s+1$ of these and one of them is the chosen $a_i$. Denote these as $x_i^1, x_i^2, \ldots, x_i^{2s+1}, \ldots$. Consider the $k$-tuple $a_1a_2\cdots a_{i-1} x_i^j a_{i+1}\cdots a_k$ and its possible placement in the empty $s\times(s+1)$ subarray formed by the intersection of the $s$ rows through $A_i$ and $s+1$ columns through $D$. Such a placement must be impossible since $M$ is maximal, so each $x_i^j$ must lie in the rows through $A_i$ or in the columns of $D$.

Suppose that a particular $x_i^j$ lies in $\gamma$ rows through $A_i$ and in $\delta$ columns of $D$. These entries will cover \[\gamma(s+1)+\delta s-\gamma\delta = s(s+1)-(s-\gamma)(s+1-\delta)\] of the $s(s+1)$ empty cells. So either $\gamma=s$ or $\delta=s+1$.

If $f(x_i^j)=s$ then $\gamma+\delta\le s$, and so $\gamma=s$ and $\delta=0$. It follows that each $x_i^j$ with frequency $s$ occurs in all the rows through $A_i$. The same argument applies to columns through $A_i$, so such points must occur entirely within the subsquare $A_i$ and they cannot occur in $D$. There is at least one such point, namely $a_i$. It also follows that $D$ can only contain points with frequency $s+1$.

Now consider an $i^{\mathrm{th}}$ coordinate entry $x_i^j$ that occurs in $D$. If it lies in $s$ or fewer columns of $D$ then it must also lie in all $s$ rows through $A_i$ and so it must occur in only a single column of $D$. The same argument applies to rows of $D$ so if $x_i^j$ occurs only once in $D$, then all $s$ other occurrences lie in $A_i$. The only other alternative is that $x_i^j$ occurs in all $s+1$ columns of $D$, and consequently in all $s+1$ rows of $D$. There can be at most $s-1$ $x_i^j$\,s that only occur once in $D$ because their $s(s-1)$ entries in $A_i$ along with the $s$ entries $a_i$ would fill $A_i$.

If there are $\xi$ entries $x_i^j$ occurring only once in $D$, and $\eta$ entries occurring $s+1$ times in $D$ then $\xi+(s+1)\eta=(s+1)^2$. But $\xi\le s-1$, so $\eta\ge \lceil s+1-(s-1)/(s+1)\rceil$, giving $\eta \ge s+1$. Hence all cells in $D$ are filled by $(s+1)$ distinct $x_i^j$\,s, which must form a Latin square in the $i^{\mathrm{th}}$ coordinates in region $D$.

\medskip

\noindent\textbf{Case (b): $A_i$ has dimension $(s+1)\times (s+1)$.} Consider the $i^{\mathrm{th}}$ coordinate entries that do not occur in a block with any $a_j$ for $j\ne i$. There must be at least $2s+2$ of these and one of them is the chosen $a_i$. Denote these as $x_i^1, x_i^2, \ldots, x_i^{2s+2}, \ldots$. Consider the $k$-tuple $a_1a_2\cdots a_{i-1} x_i^j a_{i+1}\cdots a_k$ and its possible placement in the empty $(s+1)\times(s+1)$ subarray formed by the intersection of the $s+1$ rows through $A_i$ and $s+1$ columns through $D$. Such a placement must be impossible since $M$ is maximal, so each $x_i^j$ must lie in the rows through $A_i$ or in the columns of $D$.

Suppose that a particular $x_i^j$ lies in $\gamma$ rows through $A_i$ and in $\delta$ columns of $D$. These entries will cover \[\gamma(s+1)+\delta(s+1)-\gamma\delta = (s+1)^2-(s+1-\gamma)(s+1-\delta)\] of the $(s+1)^2$ empty cells. So either $\gamma=s+1$ or $\delta=s+1$. So each $x_i^j$ either lies in $s+1$ rows through $A_i$ and consequently in no columns of $D$, or it lies in $s+1$ columns of $D$ and in no rows through $A_i$. It also follows that each $x_i^j$ has frequency $s+1$ and that the number of the $x_i^j$\,s is at most (and therefore equal to) $2s+2$. If there are $\xi$ $x_i^j$\,s with entries in rows through $A_i$ and $\eta$ $x_i^j$s with entries in $D$, then $\xi+\eta=2s+2, ~\xi\le s+1$ and $\eta\le s+1$. This implies $\xi=\eta=s+1$. So, without loss of generality $x_i^j$ lies in $s+1$ rows through $A_i$ (and not in $D$) for $1\le j\le s+1$ and $x_i^j$ lies in $s+1$ columns of $D$ (and not in rows through $A_i$) for $s+2\le j\le 2s+2$. Each of these entries in $D$ must also occur in all $s+1$ rows of $D$ and so they form a Latin square in the $i^{\mathrm{th}}$ coordinates in region $D$.

Combining cases (a) and (b), we see that $D$ forms a Latin square on each of the $k$ coordinates, irrespective of the dimensions of each $A_i$. From $M$ it inherits the orthogonality property in the sense that it contains no repeated entry pairs. Hence $D$ is a $k$-OLS($s+1$). Its entries cannot occur in $A$ since they occur $s+1$ times in $D$.

Consider again Figure \ref{fig:diag2}. Delete region $D$ and the $s+1$ entries in each of the $k$ coordinates that it contains. Also delete regions $B$ and $C$, leaving just region $A$. Put $n'=n-(s+1)$, so that $A$ is an $n'\times n'$ array on $n'$ points in each coordinate. Delete any coordinate (say the $j^{\mathrm{th}}$) to form $A'$, so that $A'$ is an $n'\times n'$ array on $n'$ points in each coordinate and now with entry $(k-1)$-tuples, and it forms a $(k-1)$-OPLS($n'$). If $A'$ were not maximal, there would be some $(k-1)$-tuple that could be placed into some empty cell in $A'$. If any one of the $s+1$ $j^{\mathrm{th}}$ coordinate entries deleted with $D$ is inserted into this $(k-1)$-tuple in the missing $j^{\mathrm{th}}$ coordinate position then the resulting $k$-tuple could be placed into $M$, a contradiction. Hence  $A'$ is a $(k-1)$-MOPLS($n'$).
The number of filled cells in $A'$ is
\[F-(s+1)^2=(k-(r-1))s^2+(r-1)(s+1)^2,\]
and all frequencies remain as $s$ or $s+1$. So by the inductive hypothesis it can be put into partitioned form. Now reinstate the deleted coordinate and regions $B,C,D$, and we have $M$ in partitioned form.

\bigskip

In order to complete the proof we must finally consider the case when $r=0$. In this case all frequencies are $s$. The argument is largely as previously except that region $D$ and each region $A_i$ has $s$ rows and $s$ columns. It is still the case that region $D$ must be full. The proof that region $D$ forms a Latin square on each coordinate follows closely case (b) with $s+1$ replaced by $s$ throughout (write $s-1$ in place of every $s$), and consequently $D$ is a \linebreak $k$-OLS($s$). The final step is as before with now $n'=n-s=ks$, so that $A'$ has $ks^2$ filled cells with all frequencies remaining as $s$.
\end{proof}

A $k$-MOPLS($n$) is said to be \emph{balanced} if all frequencies are equal. A \emph{balanced minimal} system is a balanced system with the minimum possible number of filled cells. We next show that, if $n$ is sufficiently large, any such design is equivalent to one produced by Construction 1.

\begin{theorem}\label{th:balanced} Suppose that $M$ is a balanced minimal $k$-MOPLS($n$). Then, if $n$ is sufficiently large, the common frequency is $f=n/(k+1)$ (so that $n\equiv 0 \pmod{k+1}$). Such a design is, up to permutation of rows and columns, one produced by Construction 1.
\end{theorem}
\begin{proof} Suppose that the common frequency $f<n/(k+1)$, giving $n-kf>f$. As in the proof of Theorem \ref{th:partition}, a $k$-tuple $a_1a_2\ldots a_k$ may be chosen such that every pair $a_ia_j$ ($i\ne j$) is non-occurring. By permuting rows and columns, $M$ can be partitioned  as in Figure \ref{fig:ABCDg}, with region $A$ containing all occurrences of $a_1,a_2,\ldots,a_k$ and having dimension at most $kf\times kf$. Then region $D$ has dimension greater than $f\times f$, and must be full otherwise $a_1a_2\ldots a_k$ may be inserted into an empty cell of $D$, and so rows through $D$ have frequency greater than $f$, a contradiction. Hence $f\ge n/(k+1)$.

If $n=(k+1)s+r$ where $1\le r\le k$ (so $r\ne 0$) then $f\ge s+1$ and so the number of filled cells $F\ge n(s+1)$. Now assume that $n$ is large enough to ensure that a $k$-OLS($s$) and a $k$-OLS($s+1$) both exist. Then Construction 1 produces a $k$-MOPLS($n$) with less than $n(s+1)$ filled cells. So we must have $r=0$, and then Construction 1 produces a $k$-MOPLS($n$) that is balanced with frequency $s=n/(k+1)$. So $M$, being minimal, must also have $f=s$ and consequently $n^2/(k+1)$ filled cells. Theorem \ref{th:partition} then asserts that $M$ is, up to permutation of rows and columns, one produced by Construction 1.
\end{proof}


Theorem \ref{th:partition} has the explicit assumption that there are just two possible frequencies. In certain cases, including the cases $k=2$ and $k=3$ it can be shown that this is a consequence of the other assumptions.

\begin{theorem}\label{th:maxfreq}
 Suppose that $M$ is a $k$-MOPLS($n$) with $n=(k+1)s+r$ where $0\le r\le k$, minimum frequency $m=s$, and $F=(k+1-r)s^2 +r(s+1)^2$ filled cells. If $r=0$ then all frequencies are $m$, if $r=1$ or $2$, the maximum frequency is $m+1$, and in the special case $(k,r)=(3,3)$, the maximum frequency is again $m+1$.
\end{theorem}
\begin{proof}
Choose any one of the $k+2$ groups of $M$ and for $i\ge 0$ let $n_i$ denote the number of points in that group that have frequency $s+i$. Then
\[\sum_{i\ge0} n_i=n=(k+1)s+r, \mbox{~and~} \sum_{i\ge0} (s+i)n_i=F=(k+1-r)s^2 +r(s+1)^2.\]
It follows that
\[\sum_{i\ge0} in_i=F-ns=r(s+1), \mbox{~and so~}\sum_{i\ge1} n_i\le r(s+1). \]
This gives $n_0\ge n-r(s+1)=(k+1-r)s$. It follows that every group of $M$ contains at least $(k+1-r)s$ points with frequency $s$. If $r=0$ this proves that every point has frequency $m=s$.

\smallskip

\noindent\textbf{The case $r=1$.}~Here every group of $M$ contains at least $ks$ points with frequency $s$. As in the proof of Theorem \ref{th:partition}, we can choose an entry $k$-tuple $a_1a_2\ldots a_k$ such that every pair $a_ia_j$ ($i\ne j$) is non-occurring and $f(a_i)=s$ for each $i$. By permuting rows and columns in $M$, we may take region $A$ in Figure \ref{fig:ABCDg} to be a $ks\times ks$ subsquare containing all occurrences of each $a_i$, leaving region $D$ as an $(s+1)\times (s+1)$ subsquare which must be full, otherwise $a_1a_2\ldots a_k$ could be inserted into an empty cell. Counting frequencies of the rows, we have $F\ge ks^2 +(s+1)^2$ with a strict inequality if any row has frequency greater than $s+1$. Since $F=ks^2 +(s+1)^2$, we deduce that the maximum frequency of any row is $(s+1)$, and the same argument applies immediately to columns. By permuting the groups of the design, it follows that no group has a point with frequency in excess of $s+1$.

\smallskip

\noindent\textbf{The case $r=2$.}~Here every group of $M$ contains at least $(k-1)s$ points with frequency $s$. Suppose that some group contains more than this number of points with frequency $s$; without loss of generality we may suppose that this is the $k^{\mathrm{th}}$ entry group $X_k$. As in the proof of Theorem \ref{th:partition}, we can choose an entry $(k-1)$-tuple $a_1a_2\ldots a_{k-1}$ such that every pair $a_ia_j$ ($i\ne j$) is non-occurring and $f(a_i)=s$ for each $i$. But it is also possible to choose $a_k\in X_k$, non-occurring with each $a_j$ for $j<k$ and with $f(a_k)=s$. Then by permuting rows and columns in $M$, we may take region $A$ in Figure \ref{fig:ABCDg} to be a $ks\times ks$ subsquare containing all occurrences of each $a_i$, leaving region $D$ as an $(s+2)\times (s+2)$ subsquare which must be full, otherwise $a_1a_2\ldots a_k$ could be inserted into an empty cell. Counting frequencies of the rows, we have $F\ge ks^2 +(s+2)^2$, a contradiction since $F=(k-1)s^2+2(s+1)^2$. Hence every group of $M$ contains exactly $(k-1)s$ points with frequency $s$, leaving $2s+2$ points, each with frequency at least $s+1$.

Counting frequencies in any group then gives $F\ge (k-1)s^2+2(s+1)^2$ with a strict inequality if any point has frequency greater than $s+1$. It follows that no point can have frequency greater than $s+1$.

\smallskip

\noindent\textbf{The case $(k,r)=(3,3)$.}~~Here $M$ is a 3-MOPLS($n$) with $n=4s+3$, minimum frequency $m=s$ and $F=s^2+3(s+1)^2=(n^2+3)/4$ filled cells. Every group of $M$ contains at least $s$ points with frequency $s$. Suppose that some group contains more than $s$ points with frequency $s$; without loss of generality we may suppose that this is the $2^{\mathrm{nd}}$ entry group $X_2$. As in the general case $r=2$ above, we can choose a non-occurring pair $a_1a_2$ with $a_2\in X_2$ such that $f(a_1)=f(a_2)=s$, so that $f(a_1)+f(a_2)=2s$. Referring to inequality (\ref{ineq:main2}) of Theorem \ref{th:3-MOPLS}, with $\mu\le 2s$, we have $n-2\mu\ge3$ and consequently $F\ge 9/4+n^2/4$, a contradiction. Hence every group of $M$ contains exactly $s$ points with frequency $s$, leaving $3s+3$ points, each with frequency at least $s+1$.

Counting frequencies in any group then gives $F\ge s^2+3(s+1)^2$ with a strict inequality if any point has frequency greater than $s+1$. Consequently no point can have frequency greater than $s+1$.
\end{proof}

It follows from Theorem \ref{th:maxfreq} and Corollary \ref{cor:main} that if $M$ is a minimal 2-MOPLS($n$) with $n=3s+r$, $0\le r\le 2$, and $n\ge 21$, then the maximum frequency of any point is $s+1$. Similarly, using Theorem \ref{th:maxfreq}, Corollary \ref{cor:main2} and Theorem \ref{th:minfreq}, when $M$ is a minimal 3-MOPLS($n$) with $n=4s+r$, $0\le r\le 3$, and $n\ge 44$, then the maximum frequency of any point is $s+1$. The following corollaries are immediate consequences of these results and Theorem \ref{th:partition}.

\begin{cor}\label{cor:structure2}
Suppose that $M$ is a minimal $2$-MOPLS($n$) with $n\ge 21$. Then $M$ can be partitioned into three $2$-OLSs on different point sets as in Construction 1. The orders of the three $2$-OLSs are either $s$ or $s+1$, where $s=\lfloor n/3 \rfloor$, the three orders sum to $n$, and the number of filled cells is $\lceil n^2/3\rceil$.
\end{cor}

\begin{cor}\label{cor:structure3}
Suppose that $M$ is a minimal $3$-MOPLS($n$) with $n\ge 44$. Then $M$ can be partitioned into four $3$-OLSs on different point sets as in Construction 1. The orders of the four $3$-OLSs are either $s$ or $s+1$, where $s=\lfloor n/4 \rfloor$, the four orders sum to $n$, and the number of filled cells is $\lceil n^2/4 \rceil$ for $n\equiv 0,1,3 \pmod{4}$ and $n^2/4 +1$ for $n\equiv 2 \pmod{4}$.
\end{cor}

\section{Concluding remarks}
Corollaries \ref{cor:structure2} and \ref{cor:structure3} establish the unique structures of minimal 2-MOPLS($n$) and 3-MOPLS($n$) for sufficiently large $n$.  Theorem \ref{th:balanced} does the same for balanced minimal $k$-MOPLS($n$) for all $k$. We conjecture that, for sufficiently large $n$, any minimal $k$-MOPLS($n$) must have the structure of Construction 1, which would consequently determine the minimum number of filled cells. In graph-theoretic terms this would imply that, for sufficiently large $n$, the minimum number of $K_r$ subgraphs that must be removed from  $K_{r\times n}$ in order to leave a $K_r$-free graph is approximately $n^2/(r-1)$, with the precise number depending on the residue of $n$ modulo $r-1$.

Some of our results can be extended to a more general graph-theoretic setting than that corresponding to orthogonal partial Latin squares. Let $G(r,n)$ be an $r$-partite graph with parts $X_i$ for $i=1,2,\ldots,r$, all of size $n$. Define $f(x,X_i)$ to be the number of edges between a vertex $x$ and a part $X_i$. We may define $G(r,n)$ to be \emph{vertex-balanced} if, for each vertex $x$, $f(x,X_i)=f(x,X_j)$ whenever $x\in X_k$ and $k\ne i,j$. If $G(r,n)$ is vertex-balanced then the number of edges $F$ joining parts $X_i$ and $X_j$ ($i\ne j$) is independent of the choice of $i$ and $j$. When $G(k+2,n)$ corresponds
to a $k$-OPLS, then $F$ is the number of filled cells. Our proof of Theorem \ref{th:main} can be adapted to show that if $G=G(4,n)$ is vertex-balanced and the complementary graph $G'=K_{4\times n}-G$ contains no $K_4$, then $F\ge n^2/3$ and equality is only possible if $f(x,X_i)=n/3$ for $x\not\in X_i$.

Our final comment is that by taking the block design viewpoint, each $k$-MOPLS($n$) gives rise to an $n$-ary code of length $k+2$ by regarding the blocks as codewords. There is not room here to provide details, but for $n>3$ a 2-MOPLS($n$) is equivalent to an $n$-ary code of length 4 having minimum distance 3 and covering radius 2. And it follows from our Theorem \ref{th:main} and Corollary \ref{cor:main} that for $n\ge 21$ the minimum number of codewords in such a code is $\lceil n^2/3 \rceil$. Apart from trivial cases, the codes produced by Construction~1 for $k$-MOPLS($n$) have minimum distance $k+1$ and covering radius $k$. These codes are non-extendable in the sense that adding any $(k+2)$-tuple as a codeword reduces the minimum distance.

\end{document}